# On the structure of some modules over generalized soluble groups

L.A. Kurdachenko, I.Ya. Subbotin and
V.A. Chepurdya

**Abstract**. Let R be a ring and G a group. An R-module A is said to be ***artinian-by-(finite rank)***, if $\mathbf{Tor}_R(A)$ is artinian and $A/\mathbf{Tor}_R(A)$ has finite R-rank. The authors study $\mathbb{Z}G$-modules A such that $A/C_A(H)$ is artinian-by-(finite rank) (as a $\mathbb{Z}$-module) for every proper subgroup H.

The modules over group rings RG are classic objects of study with well established links to various areas of algebra. The case when G is a finite group has been studying in sufficient details for a long time. For the case when G is an infinite group, the situation is different. Investigation of modules over polycyclic-by-finite groups was initiated in classical works of P. Hall [**HP1954**, **HP1959**]. Nowadays the theory of modules over polycyclic - by - finite groups is highly developed and rich on interesting results. This was largely due to the fact, that a group ring RG of a polycyclic - by - finite group G over a noetherian ring R is also noetherian. The group rings over some other groups (even well-studied as for instance, the Chernikov groups), do not always have such good properties as to be noetherian. Therefore, their study requires some different approaches and restrictions. Very popular such restrictions are the classical finiteness conditions. The very first restrictions here were those who came from ring theory, namely the conditions like "to be noetherian" and "to be artinian". Noetherian and artinian modules over group rings are also very well investigated. Many aspects of the theory of artinian modules over group rings are reflected in the book [**KOS2007**]. Recently the so-called finitary approach is under intensive development. This is mainly due to the progress which its application has in the theory of infinite dimensional linear groups.

Let R be a ring, G a group and A an RG-module. For a subgroup H of G we consider the R –submodule $C_A(H)$. Then H really acts on $A/C_A(H)$. The R-factor-module $A/C_A(H)$ is called the ***cocentralizer of H in A***. The factor-group $H/C_H(A/C_A(H))$ is isomorphic to a subgroup of automorphisms group of an R-module $A/C_A(H)$. If x is an element of $C_H(A/C_A(H))$, then x acts trivially on factors of the series $\langle 0\rangle \le C_A(H) \le A$. It follows that $C_H(A/C_A(H))$ is abelian. This shows that the structure of H to a greater extent is defined by the structure of $C_H(A/C_A(H))$, and hence by the structure of the automorphisms group of the R-module $A/C_A(H)$. Let $\mathfrak{M}$ be a class of R-modules. We say that A is ***$\mathfrak{M}$-finitary module over RG***, if $A/C_A(x) \in \mathfrak{M}$ for each element $x \in G$. If R is a field, $C_G(A) = \langle 1\rangle$ and $\mathfrak{M}$ is a class of all finite dimensional vector spaces over R, then we come to the finitary linear groups. The theory of finitary linear groups is quite well developed (see, the survey [**PR1995**]).

B.A.F. Wehrfritz began to consider the cases when $\mathfrak{M}$ is the class of finite R-modules [**WB2002[1], WB2002[3], WB2002[4], WB2004[1]**], when $\mathfrak{M}$ is the class of noetherian R-modules [**WB2002[2]**], when $\mathfrak{M}$ is the class of artinian R-modules [**WB2002[4], WB2003, WB2004[1], WB2004[2], WB2005**]. The artinian-finitary modules have been considered also in the paper [**KSC2007**]. The artinian and noetherian modules can be united into the following type of modules. An R-module A is said to be ***minimax***, if A has a finite series of submodules, whose factors are either noetherian or artinian. It is not hard to show that if R is an integral domain, then every minimax R-module A includes a noetherian submodule B such that A/B is artinian. The first natural case here is the case when R = $\mathbb{Z}$ is the ring of all integers. This case has very important applications in generalized soluble groups. Every $\mathbb{Z}$-minimax module M has he following important property: $\mathbf{r}_{\mathbb{Z}}(M)$ is finite and **Tor**(M) is an artinian $\mathbb{Z}$-module.

Let R be an integral domain and A be an R-module. An analogue of the concept of the dimension for modules over integral domains is the concept of R - rank. One of the essential differences of R - modules and vector spaces is that some elements of A can have a non-zero annihilator in the module. Put $\mathbf{Tor}_R(A) = \{ a \in A \mid \mathbf{Ann}_R(a) \neq \langle 0 \rangle \}$. It is not hard to see that $\mathbf{Tor}_R(A)$ is an R-submodule of A. In reality, the concept of R-rank works only for the factor-module $A/\mathbf{Tor}_R(A)$. In particular, the finiteness of R - rank does not affect the submodule $\mathbf{Tor}_R(A)$. We say that an R-module A is an ***artinian-by-(finite rank)***, if $\mathbf{Tor}_R(A)$ is artinian and $A/\mathbf{Tor}_R(A)$ has finite R-rank. In particular, if an artinian-by-(finite rank) module A is R-torsion-free, then it could be embedded into a finite dimensional vector space (over the field of fractions of R). If A is R-periodic, then it is artinian.

Let G be a group, A an RG-module, and $\mathfrak{M}$ a class of R-modules. Put

$$\mathcal{C}_{\mathfrak{M}}(G) = \{ H \mid H \text{ is a subgroup of } G \text{ such that } A/C_A(H) \in \mathfrak{M} \}.$$

If A is an $\mathfrak{M}$-finitary module, then $\mathcal{C}_{\mathfrak{M}}(G)$ contains every cyclic subgroup (moreover, every finitely generated subgroup whenever $\mathfrak{M}$ satisfies some natural restrictions). It is clear that the structure of G depends significantly on which important subfamilies of the family $\mathcal{L}(G)$ of all proper subgroups of G include $\mathcal{C}_{\mathfrak{M}}(G)$. The first natural question that arises here is the following: What is the structure of a group G in which $\mathcal{L}(G) = \mathcal{C}_{\mathfrak{M}}(G)$ (in other words, the cocentralizer of every proper subgroup of G belongs to $\mathfrak{M}$)? In the current article, we consider the case when R = $\mathbb{Z}$ and $\mathfrak{M}$ is the class of all artinian-by-(finite rank) modules. This examination is conducted for the groups belonging to the very wide class of the locally generalized radical groups.

A group G is called ***generalized radical***, if G has an ascending series whose factors are locally nilpotent or locally finite. Hence a generalized radical group G either has an ascendant locally nilpotent subgroup or an ascendant locally finite subgroup. In the first case, the locally nilpotent radical ***Lnr(G)*** of G is non-identity. In the second case, it is not hard to see that G includes a non-identity normal locally finite subgroup. Clearly in every group G the subgroup ***Lfr(G)*** generated by all normal locally finite subgroups is the largest normal locally finite subgroup (the ***locally finite radical***). Thus every generalized radical group has an ascending series of normal subgroups with locally nilpotent or locally finite factors.

Observe also that a periodic generalized radical group is locally finite, and hence a periodic locally generalized radical group is also locally finite.

**THEOREM.** *Let G be a locally generalized radical group and A a* $\mathbb{Z}G$*-module. If the factor-module* $A/C_A(H)$ *is artinian-by-(finite rank) for every proper subgroup H of G, then either* $A/C_A(G)$ *is artinian-by-(finite rank) or* $G/C_G(A)$ *is a cyclic or quasicyclic p-group for some prime p.*

**COROLLARY M.** *Let G be a locally generalized radical group and A a* $\mathbb{Z}G$*-module. If a factor-module* $A/C_A(H)$ *is minimax for every proper subgroup H of G, then either* $A/C_A(G)$ *is minimax or* $G/C_G(A)$ *is a cyclic or quasicyclic p-group for some prime p.*

**COROLLARY N.** *Let G be a locally generalized radical group and A a* $\mathbb{Z}G$*-module. If a factor-module* $A/C_A(H)$ *is finitely generated for every proper subgroup H of G, then either* $A/C_A(G)$ *is finitely generated or* $G/C_G(A)$ *is a cyclic or quasicyclic p-group for some prime p.*

**COROLLARY A.** *Let G be a locally generalized radical group and A a* $\mathbb{Z}G$*-module. If a factor-module* $A/C_A(H)$ *is artinian for every proper subgroup H of G, then either* $A/C_A(G)$ *is artinian or* $G/C_G(A)$ *is a cyclic or quasicyclic p-group for some prime p.*

We also will show that for every quasicyclic group one can find a $\mathbb{Z}G$-module A such that $C_G(A) = \langle 1 \rangle$ but the factor-module $A/C_A(H)$ is artinian-by-(finite rank) for every proper subgroup H of G.

Let R be a ring and $\mathfrak{M}$ a class of R-modules. Then $\mathfrak{M}$ is said to be a ***formation*** if it satisfies the following conditions:

(**F1**) if $A \in \mathfrak{M}$ and B is an R-submodule of A, then $A/B \in \mathfrak{M}$;

(**F2**) if A is an R-module and $B_1, \dots, B_k$ are R-submodules of A such that $A/B_j \in \mathfrak{M}$, $1 \le j \le k$, then $A/(B_1 \cap \dots \cap B_k) \in \mathfrak{M}$.

**LEMMA 1.** *Let R be a ring,* $\mathfrak{M}$ *a formation of R-modules, G a group and A an RG-module.*

*(i) If L, H are subgroups of G such that* $L \le H$ *and* $A/C_A(H) \in \mathfrak{M}$*, then* $A/C_A(L) \in \mathfrak{M}$.

*(ii) If L, H are subgroups of G whose cocentralizers belong to* $\mathfrak{M}$*, then* $A/C_A(\langle H, L \rangle) \in \mathfrak{M}$.

**Proof.** The inclusion $L \le H$ implies the inclusion $C_A(L) \ge C_A(H)$. Since $A/C_A(H) \in \mathfrak{M}$ and M is a formation, $A/C_A(L) \in \mathfrak{M}$. Clearly $C_A(\langle H, L \rangle) \le C_A(H) \cap C_A(L)$. Since M is a formation, $A/(C_A(H) \cap C_A(L)) \in \mathfrak{M}$. Then and $A/C_A(\langle H, L \rangle) \in \mathfrak{M}$.

Clearly the class of artinian-by-(finite rank) modules over an integral domain R is a formation. So we obtain

**COROLLARY.** *Let R be a ring, G a group and A an RG-module.*

*(i) If $L, H$ are subgroups of $G$ such that $L \le H$ and the factor-module $A/C_A(H)$ is artinian-by-(finite rank), then $A/C_A(L)$ is also an artinian-by-(finite rank) module.*

*(ii) If $L, H$ are subgroups of $G$, whose cocentralizers are artinian-by-(finite rank) modules, then $A/C_A(< H, L >)$ is also artinian-by-(finite rank).*

A group $G$ is said to be ***𝔉-perfect*** if $G$ does not include proper subgroups of finite index.

**LEMMA 2.** *Let $G$ be a locally generalized radical group and $A$ be a $\mathbb{Z}G$-module. Suppose that $A$ includes a $\mathbb{Z}G$-submodule $B$ which is artinian-by-(finite rank). Then the following assertions hold:*

*(i) $G/C_G(B)$ is soluble-by-finite.*

*(ii) If $G/C_G(B)$ is periodic, then it is nilpotent-by-finite.*

*(iii) If $G/C_G(B)$ is 𝔉-perfect and periodic, then it is abelian. Moreover $[[B, G], G] = <0>$.*

**Proof.** Without loss of generality we can suppose that $C_G(B) = <1>$. Being artinian-by-(finite rank), $B$ has a series of $G$-invariant subgroups $<0> \le D \le K \le B$ where $D$ is a divisible Chernikov subgroup, $K/D$ is finite and $B/K$ is torsion-free and has finite $\mathbb{Z}$-rank. More exactly, $D = D_1 \oplus \ldots \oplus D_n$ where $D_j$ is a Sylow $p_j$-subgroup of $D$, $1 \le j \le n$. Clearly every subgroup $D_j$ is $G$-invariant, $1 \le j \le n$. Let $q = p_j$. The factor-group $G/C_G(D_j)$ is isomorphic to a subgroup of $\mathbf{GL}_m(\mathbb{Z}_{q^\infty})$ where $\mathbb{Z}_{q^\infty}$ is the ring of integer $q$-adic numbers and $m$ satisfies $q^m = |\Omega_1(D_j)|$. Let $F$ be a field of fractions of $\mathbb{Z}_{q^\infty}$, then $G/C_G(D_j)$ is isomorphic to a subgroup of $\mathbf{GL}_m(F)$. Note that $\mathbf{char}(F) = 0$. Being locally generalized radical, $G/C_G(D_j)$ does not include the non-cyclic free subgroup; thus an application of Tits Theorem (see, for example, [**WB1973**, Corollary 10.17]) shows that $G/C_G(D_j)$ is soluble-by-finite. If $G$ is periodic, then $G/C_G(D_j)$ is finite (see, for example, [**WB1973**, Theorem 9.33]). It is valid for each $j$, $1 \le j \le n$. We have $C_G(D) = \bigcap_{1 \le j \le n} C_G(D_j)$, therefore using Remak's Theorem we obtain the imbedding

$$G/C_G(D) \hookrightarrow \mathbf{Dr}_{1 \le j \le n}\, G/C_G(D_j),$$

which shows that $G/C_G(D)$ is also soluble-by-finite (respectively finite). Since $K/D$ is finite, $G/C_G(K/D)$ is finite. Finally, $G/C_G(B/K)$ is isomorphic to a subgroup of $\mathbf{GL}_r(\mathbb{Q})$, where $\mathbf{r} = \mathbf{r}_{\mathbf{Z}}(B/K)$. Using again the fact that $G/C_G(A/K)$ does not include the non-cyclic free subgroup and Tits Theorem (respectively Theorem 9.33 of the book [**WB1973**]), we obtain that $G/C_G(B/K)$ is soluble-by-finite (respectively finite). Put

$$Z = C_G(D) \cap C_G(K/D) \cap C_G(B/K).$$

Then $G/Z$ is embedded in $G/C_G(D) \times G/C_G(K/D) \times G/C_G(B/K)$, in particular, $G/Z$ is soluble-by-finite (respectively finite). If $x \in Z$, then $x$ acts trivially in every factors of a series $<0> \le D \le K \le A$. By Kaluzhnin's Theorem [**KL1953**] $Z$ is nilpotent. It follows that $G$ is soluble-by-finite (respectively nilpotent-by-finite).

Suppose now that $G$ is an 𝔉-perfect group. Again consider the series of

G-invariant subgroups $<0> \le K \le B$. Being abelian and Chernikov, K is an union of ascending series

$$<0> = K_0 \le K_1 \le \ldots \le K_n \le K_{n+1} \le \ldots$$

of G-invariant finite subgroups $K_n$, $n \in \mathbf{N}$. Then the factor-group $G/C_G(K_n)$ is finite, $n \in \mathbf{N}$. Since G is $\mathfrak{F}$-perfect, $G = C_G(K_n)$ for each $n \in \mathbf{N}$. The equation

$K = \cup_{n \in \mathbf{N}} K_n$ implies that $G = C_G(K)$. By the above $G/C_G(B/K)$ is soluble-by-finite, and being $\mathfrak{F}$-perfect, it is soluble. Then $G/C_G(B/K)$ includes normal subgroups U, V such that $C_G(B/K) \le U \le V$, $U/C_G(B/K)$ is isomorphic to a subgroup of $\mathbf{UT}_r(\mathbf{Q})$, V/U includes a free abelian subgroup of finite index [**CVS1954**, Theorem 2]. Since $G/C_G(B/K)$ is $\mathfrak{F}$-perfect, it follows that $G/C_G(B/K)$ is torsion-free. Being periodic, $G/C_G(B/K)$ must be identity. In other words, $G = C_G(B/K)$. Hence G acts trivially in every factors of the series $<0> \le K \le A$, so that $[[B, G], G] = <0>$ and using again Kaluzhnin's Theorem [**KL1953**], we obtain that G is abelian.

**COROLLARY.** *Let G be a group and A a* $\mathbf{Z}$*G-module. If the factor-module* $A/C_A(G)$ *is artinian-by-(finite rank), then every locally generalized radical subgroup of* $G/C_G(A)$ *is soluble-by-finite, and every periodic subgroup of* $G/C_G(A)$ *is nilpotent-by-finite.*

Indeed, **Lemma 2** shows that $G/C_G(A/C_A(G))$ is soluble-by-finite. Every element $x \in C_G(A/C_A(G))$ acts trivially in the factors of the series $<0> \le C_A(G)) \le A$. It follows that $C_G(A/C_A(G))$ is abelian. Suppose now that $H/C_G(A)$ is a periodic subgroup. Since $A/C_A(G)$ is an $\mathbf{A}_3$-module, A has a series of H-invariant subgroups $<0> \le C_A(G) \le D \le K \le A$ where $D/C_A(G)$ is a divisible Chernikov subgroup, K/D is finite and A/K is torsion-free and has finite $\mathbf{Z}$-rank. In **Lemma 2** we have already proved that $G/C_G(D/C_A(G))$, $G/C_G(K/D)$ and $G/C_G(A/K)$ are finite. Let $Z = C_G(D/C_A(G)) \cap C_G(K/D) \cap C_G(A/K)$. Then G/Z is finite. If $x \in Z$, then x acts trivially in every factors of the series $<0> \le C_A(G) \le D \le K \le A$. By Kaluzhnin's Theorem [**KL1953**] Z is nilpotent.

Next result is well-known, but we were not able to find an appropriate reference, so we prove it here.

**LEMMA 3.** *Let G be an abelian group. Suppose that* $G \ne KL$ *for arbitrary proper subgroups K, L. Then G is a cyclic or quasicyclic p-group for some prime p.*

**Proof.** If G is finite, then it is not hard to see that G is a cyclic p-group for some prime p. Therefore suppose that G is infinite. If G is periodic, then obviously G is a p-group for some prime p. Let B be a basic subgroup of G, that is B is a pure subgroup of G such that B is a direct product of cyclic p-subgroups and G/B is divisible. The existence of such subgroups follows from [**FL1970**, Theorem 32.3]. Since G/B is divisible, $G/B = \mathbf{Dr}_{\lambda \in \Lambda} D_\lambda$ where $D_\lambda$ is a quasicyclic subgroup for every $\lambda \in \Lambda$ (see, for example, [**FL1970**, Theorem 23.1]). Our condition shows that G/B is a quasicyclic group. In particular, if $B = <1>$, that G is a quasicyclic group. Assume that $B \ne <1>$. If B is a bounded subgroup, then $G = B \times C$ for some subgroup C (see, for example, [**FL1970**, Theorem 27.5]), and we obtain a contradiction. Suppose that B is not bounded. Then B includes a subgroup $C = \mathbf{Dr}_{n \in \mathbf{N}} < c_n >$ such that $B = C \times U$ for some subgroup U and $| c_n | = p^n$, $n \in \mathbf{N}$. Let $E = < c_n^{-1} c_{n+1}^p \mid n \in \mathbf{N} >$. Then the factor-group C/E is quasicyclic, so that B/EU is also quasicyclic. It follows that G/EU is a direct product of two quasicyclic subgroups, and again come to contradiction. This contradiction shows that B = <1>, which proves a result.

**COROLLARY.** *Let G be a soluble group. Suppose that G is not finitely generated and $G \neq < K, L >$ for arbitrary proper subgroups K, L. Then G/[G, G] is a quasicyclic p-group for some prime p.*

If G is a group, then by **Tor**(G) we will denote the maximal normal periodic subgroup of G. We recall that if G is a locally nilpotent group, then **Tor**(G) is a (characteristic) subgroup of G and G/**Tor**(G) is torsion-free.

## PROOF OF THE MAIN THEOREM.

Again suppose that $C_G(A) = <1>$. Suppose that G is a finitely generated group. Then we can choose a finite subset M such that $G = < M >$, but $G \neq < S >$ for every subset $S \neq M$. If $|M| > 1$, then $M = \{x\} \cup S$ where $x \notin S$ and $S \neq \varnothing$. It follows that $< S > = U \neq G$, thus $A/C_A(U)$ is artinian-by-(finite rank). The factor $A/C_A(x)$ is also artinian-by-(finite rank), and **Corollary to Lemma 1** shows that $< x, U > = < x, S > = G$ has an artinian-by-(finite rank) cocentralizer.

Suppose that $M = \{ y \}$, that is $G = < y >$ is a cyclic group. If y has infinite order, then $< y > = < y^p >< y^q >$ where p, q are primes, $p \neq q$, and **Corollary to Lemma 1** again implies that $A/C_A(G)$ is artinian-by-(finite rank). Finally, if y has finite order, but |y| is a not a power of some prime, then < y > is a product of two proper subgroups, and **Corollary to Lemma 1** implies that $A/C_A(G)$ is artinian-by-(finite rank).

Assume now that G is not finitely generated and $A/C_A(G)$ is not artinian-by-(finite rank). Suppose that G includes a proper subgroup of finite index. Then G includes a proper normal subgroup H of finite index. We can choose a finitely generated subgroup F such that G = HF. Since G is not finitely generated, $F \neq G$. It follows that cocentralizers of both subgroups H and F are artinian-by-(finite rank). **Corollary to Lemma 1** shows that FH = G has an artinian-by-(finite rank) cocentralizer, and we obtain a contradiction. This contradiction shows that G is an $\mathfrak{F}$-perfect group.

If H is a proper subgroup of G, then **Corollary to Lemma 2** shows that H is soluble-by-finite. In particular, G is locally (soluble-by-finite). By Theorem A of the paper [**DES2005**], G includes a normal locally soluble subgroup L such that G/L is finite or locally finite simple group. Since G is an $\mathfrak{F}$-perfect group, then in the first case G = L, i.e. G is locally soluble. Consider the second case. Put $C = C_A(L)$. In a natural way, we can consider C as $\mathbf{Z}$(G/L)-module. $C_{G/L}(C)$ is a normal subgroup of G/L. Since G/L is a simple group, then either $C_{G/L}(C)$ is the identity subgroup or $C_{G/L}(C) = G/L$. In the second case $C \leq C_A(G)$ and $A/C_A(G)$ is artinian-by-(finite rank). This contradiction shows that $C_{G/L}(C) = <1>$. Let H/L be an arbitrary proper subgroup of G/L. Then H is a proper subgroup of G, therefore $A/C_A(H)$ is artinian-by-(finite rank). It follows that $C/(C \cap C_A(H))$ is also artinian-by-(finite rank). Clearly $C_C(H/L) \leq C \cap C_A(H)$, so that $C/C_C(H/L)$ is artinian-by-(finite rank). Since H/L is periodic, it is nilpotent-by-finite by **Corollary to Lemma 2**. In other words, every proper subgroup of $G/L$ is nilpotent-by-finite. Using now Theorem A of the paper [**NP1997**], we obtain that either $G/L$ is soluble-by-finite or a p-group for some prime p. In any case, $G/L$ cannot be an infinite simple group. This contradiction shows that $G$ is locally soluble. Being an infinite locally soluble group, $G$ has a non-identity proper normal subgroup. **Corollary to Lemma 2** shows that this subgroup is soluble. It follows that $G$ includes a non-identity normal abelian subgroup. In turn, it follows that the locally nilpotent radical $R_1$ of $G$ is non-identity. Suppose that $G \neq R_1$. Being $\mathfrak{F}$-perfect, $G/R_1$ is infinite. Using the above arguments, we obtain that the locally nilpotent radical $R_2/R_1$ of $G/R_1$ is non-identity. If $G \neq R_2$, then the locally nilpotent radical $R_3/R_2$ of $G/R_2$ is non-identity, and so on. Using ordinary induction, we obtain that $G$ is a radical group. Suppose that the upper radical series of $G$ is infinite and consider its term $R_\omega$, where $\omega$ is a first infinite ordinal. By its choice, $R_\omega$ is not soluble. Then **Corollary to Lemma 2** shows that $R_\omega = G$.

Since $R_n$ is a proper subgroup of $G$, $A/C_A(R_n)$ is artinian-by-(finite rank), $n \in \mathbf{N}$. $R_n$ is normal in $G$, therefore $C_A(R_n)$ is a $\mathbf{Z}G$-submodule. **Lemma 3** shows that $G/C_G(A/C_A(R_n))$ is abelian. Suppose that there exists a positive integer $m$ such that $G \neq C_G(A/C_A(R_m))$, then $[G, G]$ is a proper subgroup of $G$. Application of **Corollary to Lemma 2** shows that $[G, G]$ is soluble, thus and $G$ is soluble. This contradiction proves the equality $G = C_G(A/C_A(R_n))$. In other words, $[A, G] \leq C_A(R_n)$. Since it is valid for each $n \in \mathbf{N}$, $[A, G] \leq \bigcap_{n \in \mathbf{N}} C_A(R_n)$. The equation $G = \bigcup_{n \in \mathbf{N}} R_n$ implies that $C_A(G) = \bigcap_{n \in \mathbf{N}} C_A(R_n)$. Hence $[A, G] \leq C_A(G)$. Thus $G$ acts trivially in both factors $C_A(G)$ and $A/C_A(G)$. Using again Kaluzhnin's Theorem [**KL1953**], we obtain that $G$ is abelian. Contradiction. This contradiction proves that $G$ is soluble.

Let $D = [G, G]$. Then by **Corollary to Lemma 3** $G/D$ is a quasicyclic p-group for some prime p. It follows that $G$ has an ascending series of normal subgroups

$$D = K_0 \leq K_1 \leq \ldots \leq K_n \leq K_{n+1} \leq \ldots$$

such that $K_n/D$ is a cyclic group of order $p^n$, $n \in \mathbf{N}$, and $G = = \bigcup_{n \in \mathbf{N}} K_n$. Every subgroup $K_n$ is proper and normal in $G$, therefore $C_A(K_n)$ is a $\mathbf{Z}G$-submodule and $A/C_A(K_n)$ is artinian-by-(finite rank). **Lemma 2** shows that $[[A, G], G] \leq C_A(K_n)$. It is valid for each $n \in \mathbf{N}$, and therefore $[[A, G], G] \leq \bigcap_{n \in \mathbf{N}} C_A(K_n)$. The equation $G = \bigcup_{n \in \mathbf{N}} K_n$ implies that $C_A(G) = \bigcap_{n \in \mathbf{N}} C_A(R_n)$. Hence $[[A, G], G] \leq C_A(G)$. It follows that $G$ acts trivially in factors $C_A(G)$, $[A, G]/C_A(G)$ and $A/[A, G]$. Application of Kaluzhnin's Theorem [**KL1953**] shows that $G$ is nilpotent of class at most 2.

If $G$ is abelian, then the result follows from **Lemma 3**. Suppose that $G$ is non-abelian. Let $T = \mathbf{Tor}(G)$. If we suppose that $T \neq G$, then $G/T$ is a non-identity torsion-free nilpotent group. In particular, $G/T$ has a non-identity torsion-free abelian factor-group, which contradicts **Corollary to Lemma 3**. This contradiction shows that $G$ is a periodic group. Moreover, $G$ is a p-group. Since $G$ is nilpotent of class 2, then $[G, G] \leq \zeta(G)$. In particular, $G/\zeta(G)$ is a quasicyclic group. In this case, $[G, G]$ is a Chernikov subgroup (see, for example, [**KOS2007**, Theorem 23.1]). It follows that entire group $G$ is Chernikov. Being $\mathfrak{F}$-perfect, $G$ is abelian, which completes the proof.

We will construct the following example showing that for every quasicyclic p-group G there exists a $\mathbb{Z}G$-module A such that every proper subgroup of G has artinian-by-(finite rank) cocentralizer.

Let p be a prime and G a quasicyclic p-group, that is $G = \langle g_n \mid n \in \mathbb{N} \rangle$, where $g_1^p = 1$, $g_2^p = g_1$, $g_{n+1}^p = g_n$, $n \in \mathbb{N}$. Let $A_j$, B be an additively written quasicyclic p-group, $j \in \mathbb{N}$, $A = \oplus_{j \in \mathbb{N}} A_j \oplus B$. Let $B = \langle b_n \mid n \in \mathbb{N} \rangle$, where $pb_1 = 0$, $pb_2 = b_1$, $pb_{n+1} = b_n$, $n \in \mathbb{N}$. Define the mapping $\gamma_1: A \longrightarrow A$, by the following rule:

$$\gamma_1(a) = \begin{cases} a + b_1, \text{ if } a \in A_1 \\ a, \text{ if } a \in \oplus_{j>1} A_j \oplus B. \end{cases}$$

Clearly $\gamma_1$ is an automorphism of A such that $\gamma_1^p = \varepsilon$ is an identity automorphism of A. Define the mapping $\gamma_2: A \longrightarrow A$, by the following rule:

$$\gamma_2(a) = \begin{cases} a + b_2, \text{ if } a \in A_1 \\ a + b_1, \text{ if } a \in A_2 \\ a, \text{ if } a \in \oplus_{j>2} A_j \oplus B. \end{cases}$$

Clearly $\gamma_2$ is an automorphism of A such that $\gamma_2^p = \gamma_1$. And in a similar way, if n is a positive integer, then define the mapping $\gamma_n: A \longrightarrow A$, by the following rule:

$$\gamma_n(a) = \begin{cases} a + b_n, \text{ if } a \in A_1 \\ a + b_{n-1}, \text{ if } a \in A_2, \\ \dots\dots\dots\dots\dots \\ a + b_1, \text{ if } a \in A_n, \\ a, \text{ if } a \in \oplus_{j>n} A_j \oplus B. \end{cases}$$

It is not hard to prove that every such mapping is an automorphism of A such that $\gamma_1^p = \varepsilon$, $\gamma_2^p = \gamma_1$, $\gamma_{n+1}^p = \gamma_n$, $n \in \mathbb{N}$. In other words, $\Gamma = \langle \gamma_n \mid n \in \mathbb{N} \rangle$ is a quasicyclic p-group. Let f be an isomorphism of G on $\Gamma$ such that $f(g_n) = \gamma_n$, $n \in \mathbb{N}$. Define the action of G on A by the rule $ag_n = \gamma_n(a)$, $n \in \mathbb{N}$. In a natural way, A becomes a $\mathbb{Z}G$-module. Furthermore, this module is periodic, $C_A(g_n) = \oplus_{j>n} A_j \oplus B$, so that $A/C_A(g_n) \cong A_1 \oplus \dots A_n$, in particular, $A/C_A(g_n)$ is an artinian $\mathbb{Z}$-module, $n \in \mathbb{N}$. At the same time, $C_A(G) = B$, and therefore $A/C_A(G)$ is not artinian.

Authors' addresses:
Leonid A.Kurdachenko, Department of Algebra and Geometry,
Dnepropetrovsk National University, Gagarin Prospect 72, Dnepropetrovsk 10, Ukraine 4910, e-mail: lkurdachenko@i.ua;

Igor Ya. Subbotin, Department of Mathematics and Natural Sciences, National University, 5245 Pacific Concourse Drive, Los Angeles, California 90045, USA, e-mail: isubboti@nu.edu;

Vasiliy Chupordya, Department of Algebra and Geometry, Dnepropetrovsk National University, Gagarin Prospect 72, Dnepropetrovsk 10, Ukraine 4910, e-mail: vchepurdya@i.ua
.